\input amstex
\documentstyle{amsppt}
\refstyle{A}
\input xypic

\magnification=1200

\topmatter

\title
       Hochschild homology of iterate skew polynomial rings
\endtitle

\author
       Jorge A. Guccione and Juan J. Guccione
\endauthor

\address
         Departamento de Matem\'atica, Facultad de Ciencias Exactas y
         Naturales, Pabell\'on 1 - Ciudad Universitaria, (1428) Buenos
         Aires, Argentina. E-mail: vander \@ mate.dm.uba.ar
\endaddress

\abstract

We study the Hochschild homology of the iterated skew polynomial
rings introduced by D. Jordan in ``A simple localization of the
quantized Weyl algebra''. First, we obtain a complex, smaller that
the canonical one of Hochschild, given the homology of such an
algebra, and then, we study this complex in order to compute the
homology of some families of algebras. In particular we compute
the homology of some quantum groups, in the generic case.
\endabstract

\subjclass
          16E40, 18G99
\endsubjclass

\thanks
\endthanks

\endtopmatter

\document

\def \sub{\subseteq}

\def \ot{\otimes}
\def \ov{\overline}
\def \wt{\widetilde}
\def \ni{\noindent}
\def \ba{\bold a}
\def \bs{\bold s}

\def \bQ{\bold Q}

\def \Z{\Bbb Z}

\def \al{\alpha}
\def \be{\beta}
\def \de{\delta}
\def \ga{\gamma}
\def \la{\lambda}

\def \fsl{\frak sl}
\def \fu{\frak u}
\def \fU{\frak U}
\def \fV{\frak V}
\def \fW{\frak W}

\def \coker {\operatorname{coker}}

\def \dg {\operatorname{dg}}
\def \Gr {\operatorname{Gr}}
\def \H {\operatorname{H}}
\def \HH {\operatorname{HH}}
\def \Hom {\operatorname{Hom}}
\def \id {\operatorname{id}}
\def \Ima {\operatorname{Im}}
\def \T {\operatorname{T}}
\def \Tor {\operatorname{Tor}}
\def \U {\operatorname{U}}

\head Introduction \endhead

Let $k$ be an arbitrary commutative ring, $A$ a $k$-algebra, $\fu$
an element of $A$ and $\ga$ an automorphism of $A$ such that
$\ga(\fu)=\fu$ and $\fu a = \ga(a)\fu$, for all $a\in A$. Let
$\al$ be an automorphism of $A$ commuting with $\ga$ and let $\be
= \ga\circ \al^{-1}$. Let $S$ be the skew polynomial ring
$A[x,\al]$. Let $p$ be an inversible element of $k$. Extend $\be$
to $S$ by setting $\be(x) = px$. There exists a $\be$-derivation
$\de$ of $S$ such that $\de(A) = 0$ and $\de(x) = \fu-p\al(\fu)$.
In \cite{J1} was introduced and studied the algebra $E =
E(A,\fu,\al,p) := S[y,\be,\de]$. It is easy to see that $E$ has
underlying abelian group $A[x,y]$ and it is the extension of $A$,
generated by the variables $x,y$ and the relations $xa = \al(a)x$,
$ya = \be(a)y$ and $yx = pxy + \fu-p\al(\fu)$. The case $A$
commutative, $p =1$ and $\ga = id$ was previously introduced in
\cite{J2} and it was studied under the ring theorist point of view
in \cite{J2}, \cite{J3}, \cite{J4} and \cite{J5}. As this
definition requires a commutative base ring, it cannot be
iterated. The above generalization was in part introduced in order
to repair this defect.

The main purpose of this paper is to study the Hochschild homology
of these algebras, under the hypothesis that $A$ is $k$-flat. In
Section~1 we obtain general results about the homology of these
algebras and, in Section~2, we apply these results to compute the
homology of several families of algebras. In particular we obtain
the homology of some quantum groups in the generic case. The
homology of these last algebras was also studied in \cite{G-G1}.
The results obtained for them here, improvement and complete some
results of our first paper. By example, now we compute completely
the homology of $\Cal O_{q^2}(sok^3)$ and $\Cal O_q(M(2,k))$, in
the generic case.

\medskip

\subhead Notations \endsubhead
Before begining, we fix some notations that we will use throughout
the paper.

\roster

\smallskip

\item For each $k$-algebra $B$, we put $\ov{B} = B/k$. Moreover, given
$b\in B$ we also let $b$ denote the class of $b$ in $\ov{B}$.

\smallskip

\item For each $k$-module $V$ we write $V^{\ot^r} = V\ot\cdots\ot V$
($r$ times).

\smallskip

\item Given $a_1\ot\cdots\ot a_r \in \ov{A}^{\ot^r}$ and $1\le i< j\le r$,
we write $\ba_{ij} = a_i\ot\cdots\ot a_j$. Moreover, for each map
$\phi\:\ov{A}\to \ov{A}$, we write $\phi(\ba_{ij})= \phi(a_i)\ot\cdots\ot
\phi(a_j)$.

\smallskip

\item Given $\ba = a_0\ot\cdots\ot a_r \in M\ot \ov{A}^{\ot^r}$,
an automorphism $\phi$ of $A$, an invertible element $p\in k$ and
$0\le l\le r$, we write
$$
\align
&d_l^{\phi}(\ba)=a_0\ot \cdots\ot a_{l-1}\ot\phi^{-1}(a_l)a_{l+1}
\ot a_{l+2}\ot \cdots\ot a_r\quad\text{if $l<r$},\\
& d_r^{\phi}(\ba) = \phi^{-1}(a_r)a_0\ot a_1 \ot \cdots \ot
a_{r-1},\\
& \de_{\phi}^{x,p}(\ba) = a_0x\ot \phi^{-1}(a_1)\ot\cdots\ot
\phi^{-1}(a_r) - p^{-1}xa_0\ot a_1\ot\cdots\ot a_r\\
\intertext{and}
&\de_{\phi}^{p,y}(\ba)=p^{-1}a_0y\ot \phi^{-1}(a_1)\ot\cdots\ot
\phi^{-1}(a_r) - ya_0\ot\cdots\ot a_r.
\endalign
$$
When $\phi = id$ we put $d_l$ instead of $d_l^{\phi}$, and when $p
=1$ we put $\de_{\phi}^x$ instead of $\de_{\phi}^{x,p}$ and
$\de_{\phi}^y$ instead of $\de_{\phi}^{p,y}$.

\smallskip

\item Given $\ba = a_0\ot\cdots\ot a_r \in M\ot \ov{A}^{\ot^r}$
and $b\in \ov{A}$, we write
$$
b\star\ba = \sum_{l=0}^r(-1)^l a_0\ot\cdots\ot a_l\ot b \ot
\ga^{-1}(a_{l+1})\ot \cdots\ot \ga^{-1}(a_r).
$$

\endroster

\head 1. Hochschild homology \endhead

Let $E$ be as in the introduction and let $M$ be an $E$-bimodule.
Assume that $A$ is $k$-flat. In this section we obtain a chain
complex, giving the Hochschild homology of $E$ with coefficients
in $M$, which is simpler than the canonical one of Hochschild.
Using this result, we obtain: in Remark~1.2 a decomposition
$\HH_*(E) = \bigoplus_{r\in \Z} \HH^{(r)}_*(E)$, of the Hochschild
homology of $E$; in Proposition~1.5 an spectral sequence
converging to $\HH^{(r)}_*(E)$; and finally, in Theorem~1.8, a
very simple complex giving $\HH_*(E)$, under suitable hypothesis.

\proclaim{Theorem 1.1} The Hochschild homology $\H_*(E,M)$, of $E$
with coefficients in $M$, is the homology of the double complex
$$
Y_{**}(M):=\qquad
\diagram
Y_{01}(M) \dto^{\varphi_0} & Y_{11}(M) \lto_(0.475){\partial_{11}} \dto^{\varphi_1}
& Y_{21}(M) \lto_(0.475){\partial_{21}}\dto^{\varphi_2} &\lto_(0.35){\partial_{31}}\cdots\\
Y_{00}(M) & Y_{10}(M) \lto_(0.475){\partial_{10}} & Y_{20}(M) \lto_(0.475){\partial_{20}}
&\lto_(0.35){\partial_{30}}\cdots\\
\enddiagram
$$
where
$$
\align
&Y_{n0}(M) = (M\ot \ov{A}^{\ot n})\oplus (M\ot \ov{A}^{\ot n-1})e_1,\\
&Y_{n1}(M)=(M\ot \ov{A}^{\ot n})e_2\oplus (M\ot \ov{A}^{\ot n-1})e_1e_2,\\
&\partial_{n0}(\ba)= \sum^n_{l=0}(-1)^l d_l(\ba),\\
&\partial_{n0}(\ba e_1) = \sum^{n-2}_{l=0}(-1)^l d_l(\ba)e_1 +
(-1)^{n-1} d_{n-1}^{\al}(\ba)e_1 + (-1)^{n-1}\de_{\al}^x(\ba),\\
&\partial_{n1}(\ba e_2)=\sum^{n-1}_{l=0}(-1)^l
d_l(\ba)e_2 + (-1)^n d_n^{\be}(\ba)e_2, \\
&\partial_{n1}(\ba e_1e_2) = \sum^{n-2}_{l=0}(-1)^l d_l(\ba)e_1e_2
+ (-1)^{n-1}d_{n-1}^{\ga}(\ba)e_1e_2
+ (-1)^{n-1}\de_{\al}^{x,p}(\ba)e_2, \\
&\varphi_n(\ba e_2)=(-1)^n\de_{\be}^y(\ba),\\
&\varphi_n(\ba e_1e_2) = (-1)^n \de_{\be}^{p,y}(\ba)e_1 +
(p^{-1}\fu-\al(\fu))\star \ba.
\endalign
$$

\endproclaim

\demo{Proof} Let $\wt{E}:= E^{\ot 2}$, endowed with the
$E$-bimodule structure given by $c\cdot (a\ot b)\cdot d= ad\ot
cb$. Consider the complex $Y_{**}(\wt{E})$. The enveloping algebra
$E^e$ acts on the left over each $Y_{nv}(\wt{E})$ via
$$
(a\ot b)\cdot ((c\ot d)\ot a_1\ot \cdots\ot a_{n-u}e_1^ue_2^v)=
(ac\ot db)\ot a_1\ot \cdots\ot a_{n-u}e_1^ue_2^v
\quad\,\,\,\text{($u,v\in \{0,1\}$).}
$$
It is obvious that each $Y_{nv}(\wt{E})$ is projective relative to
the family of all $E^e$-epimor\-phisms which split as $k$-module
maps. Moreover, the boundary maps of $Y_{**}(\wt{E})$ commute with
these actions and $Y_{**}(M)=M\ot_{E^e} Y_{**}(\wt{E})$. Then, to
prove the theorem will be sufficient to check that
$Y_{**}(\wt{E})$ is a resolution of $E$ as a left $E^e$-module,
since, in this case, we will have
$$
\H_*(Y_{**}(M))= \H_*(M\ot_{E^e} Y_{**}(\wt{E}))=
\Tor^{E^e\!,k}_*(M,E)= \H_*(E,M).
$$
But, by Proposition 1.2 of \cite{G-G2}, the total complex of the
double complex
$$
\diagram
E\ot E\dto^{\psi'_0} & E\ot \ov{S}\ot E \lto_(0.55){b'_{11}}\dto^{\psi'_1}
& E\ot \ov{S}^{\ot 2}\ot E \lto_(0.525){b'_{21}}\dto^{\psi'_2} & \lto_(0.3){b'_{31}}\cdots\\
E\ot E & E\ot \ov{S}\ot E \lto_(0.55){b'_{10}} & E\ot \ov{S}^{\ot 2}\ot E \lto_(0.525){b'_{20}}
& \lto_(0.3){b'_{30}}\cdots,
\enddiagram
$$
where the vertical and horizontal arrows are the $E$-bimodule maps
defined by
$$
\allowdisplaybreaks
\align
&b'_{n0}(\bs_{0,n+1}) = \sum^n_{l=0}(-1)^l \bs_{0,l-1}\ot s_ls_{l+1}\ot \bs_{l+2,n+1}, \\
\vspace{1.5\jot}
&b'_{n1}(\bs_{0,n+1}) = \sum^{n-1}_{l=0}(-1)^l \bs_{0,l-1}\ot
s_ls_{l+1} \ot \bs_{l+1,n+1} +(-1)^n \bs_{0,n-1}\ot\be^{-1}(s_n)s_{n+1}, \\
\vspace{1.5\jot}
&\psi'_n(\bs_{0,n+1}) = (-1)^n\Bigl(s_0y \ot \be^{-1}(s_1) \ot
\dots\ot \be^{-1}(s_n)\ot s_{n+1} - \bs_{0n} \ot ys_{n+1} \\
&\phantom{\psi'_n(\bs_{0,n+1}) = } -\sum^n_{j=1} \bs_{0,j-1} \ot
\delta \circ \be^{-1}(s_j) \ot \be^{-1}(s_{j+1}) \ot \dots \ot
\be^{-1}(s_n)\ot s_{n+1}\Bigr).
\endalign
$$
is an $E^e$-projective resolution of $E$. Moreover, by
Proposition~1.3 of \cite{G-G2}, we know that the maps
$$
\theta_{*v}\:Y_{*v}(\wt{E})@>>>(E\ot \ov{S}^{\ot *}\ot
E,b'_{*v})\qquad\text{(v=0,1),}
$$
given by
$$
\align
& \theta_{nv}(a_0\ot a_{n+1} \ot \ba_{1n}e_2^v) = \ba_{0,n+1},\\
& \theta_{nv}(a_0\ot a_n \ot \ba_{1,n-1}e_1e_2^v) =
\sum_{l=0}^{n-1}(-1)^{n-l-1} \ba_{0l}\ot x\ot
\al^{-1}(\ba_{l+1,n-1})\ot a_n,
\endalign
$$
are quasi-isomorphisms. To finish the proof it suffices to chek
that $\psi'_n \circ \theta_{n1} = \theta_{n0}\circ \varphi_n$,
which is easy.\qed
\enddemo

\remark{Remark 1.2} For each $r\in \Z$, we write
$$
\align
& Y^{(r)}_{n0}(E) = \bigoplus_{i,j\ge 0\atop j-i = r} Ax^iy^j\ot
\ov{A}^{\ot n}\oplus \bigoplus_{i,j\ge 0\atop j-i = r+1}(Ax^iy^j
\ot \ov{A}^{\ot n-1})e_1,\\
& Y^{(r)}_{n1}(E) = \bigoplus_{i,j\ge 0\atop j-i = r-1}(Ax^iy^j\ot
\ov{A}^{\ot n})e_2\oplus
\bigoplus_{i,j\ge 0\atop j-i = r}(Ax^iy^j\ot \ov{A}^{\ot n-1})e_1e_2,\\
\endalign
$$
It is easy to see that each $Y^{(r)}_{**}(E)$ is a subcomplex of
$Y_{**}(E)$ and that $Y_{**}(E) = \bigoplus_{r\in \Z}
Y^{(r)}_{**}(E)$. Hence,
$$
\HH_*(E) = \bigoplus_{r\in \Z} \HH^{(r)}_*(E),
$$
where $\HH^{(r)}_*(E)$ denotes the homology of $Y^{(r)}_{**}(E)$.
\endremark

\proclaim{Lemma 1.3}  The following equalities are valid in $E$
for all $a\in A$ and all $i,j\ge 0$:

\roster

\smallskip

\item $x^iy^ja = \ga^j( \al^{i-j}(a))x^iy^j$,

\smallskip

\item $ax^iy^jx= p^ja x^{i+1}y^j - a(p^j\al^{i+1}(\fu)-
\al^{i-j+1}(\fu))x^iy^{j-1}$,

\smallskip

\item $yax^iy^j = p^i\ga(\al^{-1}(a))x^iy^{j+1} - \ga(\al^{-1}(a))
(p^i\al^i(\fu)-\fu)x^{i-1}y^j$.

\endroster

\endproclaim

\demo{Proof} It is easy to check the first equality, the second
one when $j=0$, and the third one when $i=0$. Assume that $j\ge
1$. Let $z=yx-\fu = p(xy - \al(\fu))$ be the Casimir element of
$E$. It is easy to check that $zx=pxz$, $zy =p^{-1}yz$ and
$za=\ga(a)z$ for all $a\in A$. Using these facts, which were
established in \cite{J, Subsection~2.4}, the fact that that $\ga$
commutes with $\al$ and $\ga(\fu)=\fu$ and item~1), we obtain
$$
\align
ax^iy^jx &= ax^iy^{j-1}(z + \fu)\\
& = p^{j-1}ax^izy^{j-1}  + a \al^{i-j+1}(\fu)x^iy^{j-1}\\
& = p^jax^{i+1}y^j - p^jax^i\al(\fu)y^{j-1} + a\al^{i-j+1}(\fu)x^iy^{j-1}\\
& = p^jax^{i+1}y^j -p^ja\al^{i+1}(\fu)x^iy^{j-1} + a\al^{i-j+1}(\fu)x^iy^{j-1}.\\
\endalign
$$
This prove the second equality. The third one can be proved in a
similar way.\qed
\enddemo

\proclaim{Theorem 1.4} The boundary maps of $Y^{(r)}_{**}(E)$ are
given by
$$
\align
&\partial_{n0}(ax^iy^j\ot \ba)= a\ga^j\bigl(\al^{-r}(a_1)\bigr)x^iy^j\ot \ba_{2n}
+ \sum^n_{l=1}(-1)^l d_l(ax^iy^j\ot\ba),\\
\vspace{1.5\jot}
&\partial_{n0}(ax^iy^j\ot \ba'e_1)=\biggl(a\ga^{j}\bigl(\al^{-r-1}(a'_1)\bigr)
x^iy^j\ot \ba'_{2,n-1} + \sum_{l=1}^{n-2}(-1)^l d_l (ax^iy^j\ot \ba')\\
&\phantom{\partial_{n0}(ax^iy^j\ot \ba'e_1)} + (-1)^{n-1}\al^{-1}(a'_{n-1})
ax^iy^j\ot \ba'_{1,n-2} \biggr)e_1\\
&\phantom{\partial_{n0}(ax^iy^j\ot \ba'e_1)}
+(-1)^{n-1}p^jax^{i+1}y^j\ot \al^{-1}(\ba') +(-1)^n\al(a)x^{i+1}y^j\ot \ba'\\
&\phantom{\partial_{n0}(ax^iy^j\ot \ba'e_1)} +(-1)^na\bigl(p^j\al^{i+1}(\fu)
- \al^{-r}(\fu)\bigr)x^iy^{j-1}\ot \al^{-1}(\ba'),\\
\vspace{1.5\jot}
&\partial_{n1}(ax^iy^j\ot \ba e_2)= \biggl(a\ga^j\bigl(\al^{1-r}(a_1)\bigr)x^iy^j
\ot \ba_{2n} + \sum_{l=1}^{n-1}(-1)^l d_l(ax^iy^j\ot \ba)\\
&\phantom{\partial_{n1}(ax^iy^j\ot \ba e_2)} + (-1)^n\ga^{-1}\bigl(\al(a_n)\bigr) ax^iy^j
\ot \ba_{1,n-1}\biggr)e_2,\\
\vspace{1.5\jot}
&\partial_{n1}(ax^iy^j\ot \ba'e_{1}e_{2})= \biggl(a\ga^j\bigl(\al^{-r}(a'_1)\bigr)x^iy^j
\ot \ba'_{2,n-1} + \sum_{l=1}^{n-2}(-1)^l d_l(ax^iy^j\ot \ba')\\
&\phantom{\partial_{n1}(ax^iy^j\ot \ba'e_{1}e_{2}) }
+ (-1)^{n-1}\ga^{-1}(a'_{n-1}) ax^iy^j\ot \ba'_{1,n-2} \biggr)e_1e_2\\
&\phantom{\partial_{n1}(ax^iy^j\ot \ba'e_{1}e_{2})} +(-1)^{n-1}\Bigl(p^ja
x^{i+1}y^j\ot \al^{-1}(\ba')- p^{-1}\al(a)x^{i+1}y^j\ot \ba'\\
&\phantom{\partial_{n1}(ax^iy^j\ot \ba'e_{1}e_{2})}- a\bigl(p^j\al^{i+1}(\fu)-
\al^{-r+1}(\fu)\bigr)x^iy^{j-1}\ot \al^{-1}(\ba') \biggr)e_2,
\vspace{1.5\jot}
&\varphi_n(ax^iy^j\ot \ba e_{2})= (-1)^{n}\Bigl( ax^iy^{j+1}\ot
\ga^{-1}\bigl(\al(\ba)\bigr) - p^i\ga\bigl(\al^{-1}(a)\bigr)x^iy^{j+1}\ot \ba\\
&\phantom{\varphi_n(ax^iy^j\ot \ba e_{2})} +
\ga\bigl(\al^{-1}(a)\bigr)\bigl(p^i\al^i(\fu) - \fu\bigr) x^{i-1}y^j\ot \ba \Bigr),
\vspace{1.5\jot}
&\varphi_n(ax^iy^j\ot \ba'e_1e_{2})@!@!@!@! = @!@!@!@!@!
(-1)^{n-1}\Bigl(p^i\ga\bigl(\al^{-1}(a)\bigr)x^iy^{j+1}@!@!@!@!@!\ot
\ba' - p^{-1}ax^iy^{j+1}@!@!@!@!@!\ot \ga^{-1}\bigl(\al(\ba')\bigr)\\
&\phantom{\varphi_n(ax^iy^j\ot \ba'e_1e_{2})} - \ga\bigl(\al^{-1}(a)\bigr)\bigl(
p^i\al^i(\fu) - \fu\bigr)x^{i-1}y^j\ot \ba'\Bigr)e_1\\
&\phantom{\varphi_n(ax^iy^j\ot \ba'e_1e_{2})}
+\sum_{l=0}^{n-1}(-1)^{l}ax^iy^j\ot \ba'_{1l}\ot \bigl(p^{-1}\fu-\al(\fu)\bigr)\ot
\ga^{-1}(\ba'_{l+1,n-1}),
\endalign
$$
where $\ba = a_0\ot \cdots \ot a_n\in \ov{A}^{\ot n}$ and $\ba' =
a'_0\ot \cdots \ot a'_{n-1}\in \ov{A}^{\ot n-1}$ are elementary
tensors.
\endproclaim

\demo{Proof} All the formulas can be obtained from Theorem~1.1,
applying Lemma~1.3. For instance, we have
$$
\align
\partial_{n0}(ax^iy^j\ot \ba'e_1)& = \biggl(a
x^iy^ja'_1\ot \ba'_{2,n-1} + \sum_{l=1}^{n-2}(-1)^l d_l (ax^iy^j\ot \ba')\\
& + (-1)^{n-1}\al^{-1}(a'_{n-1}) ax^iy^j\ot \ba'_{1,n-2} \biggr)e_1\\
\allowdisplaybreak
& +(-1)^{n-1}ax^iy^jx\ot \al^{-1}(\ba') +(-1)^nxax^iy^j\ot \ba'\\
&=\biggl(a\ga^{-j}\bigl(\al^{-r-1}(a'_1)\bigr) x^iy^j\ot \ba'_{2,n-1} +
\sum_{l=1}^{n-2}(-1)^l d_l (ax^iy^j\ot \ba')\\
& + (-1)^{n-1}\al^{-1}(a'_{n-1}) ax^iy^j\ot \ba'_{1,n-2} \biggr)e_1\\
&+(-1)^{n-1}p^jax^{i+1}y^j\ot \al^{-1}(\ba') +(-1)^n\al(a)x^{i+1}y^j\ot \ba'\\
& +(-1)^na(p^j\al^{i+1}\bigl(\fu)-\al^{-r}(\fu)\bigr)x^iy^{j-1}\ot \al^{-1}(\ba'),
\endalign
$$
as we want.\qed
\enddemo

Given algebra maps $f,g\:A\to A$, we let $A_f^g$ denote $A$
endowed with the $A$-bimodule structure given by $a\cdot x\cdot b
: = g(a)xf(b)$. To simplify notations we write $A_f$ instead of
$A_f^{id}$ and $A^g$ instead of $A_{id}^g$.

\proclaim{Proposition 1.5} There is a convergent spectral sequence
$E^1_{nv} \Rightarrow \H_{n+v}^{(r)}(E)$, where
$$
\align
E^1_{n0} & = \bigoplus_{j\ge \max(0,r)} \H_n(A,A_{\ga^j\circ \al^{-r}}),\\
E^1_{n1} & = \bigoplus_{j\ge \max(0,r+1)} \H_n(A,A^{\al^{-1}
}_{\ga^j\circ \al^{-1-r}})\oplus \bigoplus_{j\ge \max(0,r-1)}
\H_n(A,A^{\ga^{-1}\circ \al}_{\ga^j\circ \al^{1-r}}),\\
E^1_{n2} & = \bigoplus_{j\ge \max(0,r)} \H_n(A,A^{\ga^{-1}
}_{\ga^j\circ \al^{-r}}),\\
E^1_{nv} & = 0,\qquad\text{for all $v>2$.}
\endalign
$$
\endproclaim

\demo{Proof} Let $Y_*^{(r)}(E)$ be the total complex of
$Y_{**}^{(r)}(E)$. Let us consider the filtration $0\subseteq
F^0_*\subseteq F^1_*\subseteq F^2_* = Y_*^{(r)}(E)$ of
$Y^{(r)}_*(E)$, given by
$$
\align
F^0_n & = \bigoplus_{i,j\ge 0\atop j-i = r} Ax^iy^j\ot
\ov{A}^{\ot n},\\
F^1_n & = F^0_n \oplus \bigoplus_{i,j \ge 0\atop j-i =
r+1}(Ax^iy^j \ot \ov{A}^{\ot n-1})e_1 \oplus \bigoplus_{i,j\ge
0\atop j-i = r-1}(Ax^iy^j\ot \ov{A}^{\ot n-1})e_2,\\
F^2_n & = F^1_n \oplus \bigoplus_{i,j\ge 0\atop j-i =
r}(Ax^iy^j\ot \ov{A}^{\ot n-2})e_1e_2.
\endalign
$$
It is easy to check that the spectral sequence associate with this
filtration has the required properties.\qed

\enddemo

\proclaim{Lemma 1.6} If there exists $n_0 > 0$ such that,
for all $i\ge 0$ and all $n\ge n_0$,
$$
p^{-i}(\ga^{-1}\circ \al)^{\ot n+1} - \id_{A\ot \ov{A}^{\ot n}}\:
A\ot \ov{A}^{\ot n}\to A\ot \ov{A}^{\ot n}
$$
is a bijective map, then $\HH^{(r)}_*(E)$ is the homology of the
subcomplex $Y^{(r)}_{**}(E)_{|n_0}$ of $Y^{(r)}_{**}(E)$, defined
by
$$
Y^{(r)}_{nv}(E)_{|n_0} := \cases Y^{(r)}_{nv}(E) &\quad\text{if $n\le n_0$,}\\
0 &\quad\text{if $n > n_0$.} \endcases
$$
\endproclaim

\demo{Proof} Let $n>n_0$. Let us consider the filtration
$F_{nv}^0\subseteq F_{nv}^1\subseteq F_{nv}^2\subseteq F_{nv}^3\subseteq
\cdots$ ($v=0,1$) of $Y^{(r)}_{nv}(E)$, given by
$$
\align
& F_{n0}^l:= \bigoplus\sb{ 0\le i \le l\atop i\ge -r}
Ax^iy^{r+i} \ot \ov{A}^{\ot n}\oplus \bigoplus\sb{ 0\le i \le
l-1\atop i\ge -r-1}(Ax^iy^{r+i+1}\ot \ov{A}^{\ot n-1})e_1\\
\intertext{and}
& F_{n1}^l:= \bigoplus\sb{ 0\le i \le l\atop i\ge -r+1}
(Ax^iy^{r+i-1} \ot \ov{A}^{\ot n})e_2\oplus \bigoplus\sb{ 0\le i \le
l-1\atop i\ge -r}(Ax^iy^{r+i}\ot \ov{A}^{\ot n-1})e_1e_2.
\endalign
$$
Let $\Gr(Y_{nv}^{(r)})$ be the graded module associated with this
filtration and let
$$
\wt{\varphi}_n \: \Gr(Y_{n1}^{(r)}) \to \Gr(Y_{n0}^{(r)})
$$
be the map induced by $\varphi_n\:Y^{(r)}_{n1}(E) \to
Y^{(r)}_{n0}(E)$. By the hypothesis, $\wt{\varphi}_n$ is
bijective. Then, $\varphi_n$ is also bijective. The proposition
follows immediately from this fact.\qed

\enddemo

\proclaim{Theorem 1.7} If there is a two sides ideal $I$ of $A$,
such that:

\roster

\smallskip

\item $A= k \oplus I$,

\smallskip

\item $\alpha(I)=I$ and $\gamma(I)=I$,

\smallskip

\item $p^{-i}(\gamma^{-1}\circ \alpha)^{\ot n} - \id_{I^{\ot n}}$
is a bijective map for all $i\ge 0$ and all $n\ge 1$,

\smallskip

\endroster

\noindent then $\HH^{(r)}_*(E)$ is the homology of the double
complex
$$
\diagram
{}\save[]+<0pc,-0.5pc>\Drop{\displaystyle{\bigoplus_{i,j\ge 0\atop
j-i = r-1}} k\,x^iy^je_2}\restore
{}\save[]+<0pc,-1.5pc>\xto[1,0]+<0pc,-0.6pc>^(0.5){\ov{\varphi}_0}\restore
& {}\save[]+<6pc,-0.5pc> \Drop{\displaystyle{ \bigoplus_{i,j\ge
0\atop j-i = r} k\,x^iy^je_1e_2}} \restore \save[]+<3pc,0pc>
{}\xto[0,-1]+<3pc,0pc>_(0.45){\ov{\partial}_1}\restore
{}\save[]+<5.5pc,-1.5pc>\xto[1,0]+<5.5pc,-0.6pc>^(0.5){\ov{\varphi}_0}\restore
\\
{}\save[]+<0pc,-2.5pc>\Drop{\displaystyle{\bigoplus_{i,j\ge 0\atop
j-i = r}k\,x^iy^j}}\restore & {}\save[]+<6pc,-2.5pc>
\Drop{\displaystyle{\bigoplus_{i,j\ge 0 \atop j-i =
r+1}k\,x^iy^je_1}}\save[]+<3.35pc,-2pc> {}
\xto[0,-1]+<2.25pc,-2pc>_(0.45){\ov{\partial}_0},
\enddiagram
$$
whose boundary maps are defined by
$$
\align
&\ov{\partial}_0(x^iy^je_1)= (p^j - 1)x^{i+1}y^j -
(p^j-1)\ov{\fu}x^iy^{j-1},\\
&\ov{\partial}_1(x^iy^je_1e_2)= (p^j - p^{-1})x^{i+1}y^je_2 -
(p^j-1)\ov{\fu}x^iy^{j-1}e_2,\\
&\ov{\varphi}_0(x^iy^je_2)= (1 - p^i)x^iy^{j+1} + (p^i-1)\ov{\fu}x^{i-1}y^j,\\
&\ov{\varphi}_1(x^iy^je_1e_2)= (p^i - p^{-1})x^iy^{j+1}e_1 -
(p^i-1)\ov{\fu}x^{i-1}y^je_1,
\endalign
$$
where $\ov{\fu}$ denotes the class of $\fu$ in $k\simeq A/I$.
\endproclaim

\demo{Proof} Condition~3) implies that the hypothesis of
Lemma~1.6 is satisfied, with $n_0 = 1$. Hence,
$\HH^{(r)}_*(E)= \H_*\bigl(Y^{(r)}_{**}(E)_{|1}\bigr)$. By
condition~2), $Y^{(r)}_{**}(E)_{|1}$ has a subcomplex $Z_{**}$,
defined by
$$
\alignat2
& Z_{00} = \bigoplus_{i,j\ge 0\atop j-i = r}I\,x^iy^j, &&\qquad
Z_{10} = \bigoplus_{i,j\ge 0\atop j-i = r}A\,x^iy^j\ot \ov{A}
\oplus \bigoplus_{i,j\ge 0\atop j-i = r+1}I\,x^iy^je_1,\\
& Z_{01} = \bigoplus_{i,j\ge 0\atop j-i = r-1}I\,x^iy^je_2, &&\qquad
Z_{11} = \bigoplus_{i,j\ge 0\atop j-i = r}(A\,x^iy^j\ot \ov{A})e_2
\oplus \bigoplus_{i,j\ge 0\atop j-i = r}I\,x^iy^je_1e_2.
\endalignat
$$
The argument used to prove Lemma~1.6, shows that $Z_{**}$ is
an exact complex. So, $\HH^{(r)}_*(E)$ is the homology of
$\frac{Y^{(r)}_{**}(E)_{|1}}{Z_{**}}$, which is the complex of the
statement.\qed

\enddemo

\proclaim{Corollary 1.8} Assume that $p =1$. If the hypothesis of
Theorem~1.7 are verified, then
$$
\align
& \HH_0(E) = k[x,y],\\
& \HH_1(E) = k[x,y]e_1\oplus k[x,y]e_2,\\
& \HH_2(E) = k[x,y]e_1e_2,\\
& \HH_n(E) = 0,\qquad\text{for all $n>2$.}
\endalign
$$
\endproclaim

\demo{Proof} It follows immediately from Remark~1.2 and
Theorem~1.7.\qed
\enddemo

\head 2. Same examples \endhead

In this section we use the results obtained in Section~1 in order
to compute the Hochschild homology of some families of algebras,
that appear as iterated skew polynomial rings. In all these
examples $\ga = \id$ and $p = 1$. The main results are enounced in
Theorems~2.1.1 and 2.2.8 and Corollary~2.2.9.

Let $\bQ = \{q_{ij}:1\le i,j\le v\}$ be a set of elements of $k$,
verifying $q_{ij}q_{ji} = q_{ii} = 1$ for all $1\le i,j\le l v$.
The $v$-dimensional multiparametric quantum affine space
$k_{\bQ}[t_1,\dots,t_v]$, with parameters $\bQ$, is the $k$
algebra generated by variables $t_1,\dots,t_v$ subject to the
relations $t_jt_i = q_{ij}t_it_j$.

\specialhead 2.1. Case $A = k_{\bQ}[t_1,\dots,t_v]$ an arbitrary
multiparametric quantum affine space, $\fu\in A$ and $\al(t_i) = q
t_i$, with $q\in k\setminus\{0\}$ a non root of unity
\endspecialhead
In this case $E$ is the algebra generated over $k$ by the
variables $t_1,\dots t_v,x,y$ and the relations $t_jt_i =
q_{ij}t_it_j$, $xt_i = qt_ix$, $yt_i = q^{-1}t_iy$ and $yx = xy +
\fu - \fu(qt_1,\dots,qt_v)$.

\proclaim{Theorem 2.1.1} The Hochschild homology of $E$ is given
by:
$$
\align
& \HH_0(E) = k[x,y],\\
& \HH_1(E) = k[x,y]e_1\oplus k[x,y]e_2,\\
& \HH_2(E) = k[x,y]e_1e_2,\\
& \HH_n(E) = 0,\qquad\text{for all $n>2$.}
\endalign
$$
\endproclaim

\demo{Proof} It is an immediate consequence of Corollary~1.8.\qed
\enddemo

Theorem 2.1.1 applies to the quantum algebras $\Cal
O_{q^2}(sok^3)$ and $\Cal O_q(M(2,k))$ (see \cite{J2} and
\cite{S1} for the definitions).

\specialhead 2.2. Case $A = k[t]$ with $k$ a characteristic $0$
field, $\fu\in A$ and $\al(t) = t+\la$, with $\la\in
k\setminus\{0\}$
\endspecialhead
Note that in this case $E$ is the algebra generated over $k$ by
the variables $t,x,y$ and the relations $xt = (t+\la)x$, $yt =
(t-\la)y$ and $yx = xy + \fu - \fu(t+\la)$.

It is well known that the $A$-bimodule complex $X'_*(A):= A\ot A
@<d'_1<< A\ot A$, where $d'_1(1\ot 1) = 1\ot t - t\ot 1$, is an
$A^e$-projective resolution of $A$. Moreover, there are chain
complexes maps
$$
\theta'_*\:X'_*(A)\to (A\ot\ov{A}^{\ot *}\ot A,b'_*) \qquad
\text{and} \qquad \vartheta'_*\:(A\ot\ov{A}^{\ot *}\ot A,b'_*)\to
X'_*(A),
$$
given by $\theta'_0 = \vartheta'_0 = \id_{A\ot A}$,
$\theta'_1(1\ot 1) = 1\ot t\ot 1$ and $\vartheta'_1(1\ot t^n\ot 1)
= \sum_{i=0}^{n-1} t^i\ot t^{n-i-1}$. Let $f,g\:A\to A$ be
automorphisms and let $A_f^g$ be as in Proposition~1.5. By
tensoring $A_f^g$ on the left over $A^e$ with $X'_*(A)$, we obtain
the complex
$$
X_*(A_f^g):= A_f^g @<d_1<< A_f^g,\qquad\text{where $d_1(P) =
(f(t)-g(t))P$,}
$$
which gives the Hochschild homology of $A$ with coefficients in
$A_f^g$. The maps $\theta'_*$ and $\vartheta'_*$ induce
quasi-isomorphisms
$$
\theta^{fg}_*\:X_*(A_f^g)\to (A_f^g\ot\ov{A}^{\ot *},b_*) \qquad
\text{and} \qquad \vartheta^{fg}_*\:(A_f^g\ot\ov{A}^{\ot
*},b_*)\to X_*(A_f^g).
$$
Explicitly, we have
$$
\theta^{fg}_0 = \vartheta^{fg}_0 = \id_{A_f^g},\quad
\theta^{fg}_1(P) = P\ot t\quad \text{and}\quad
\vartheta^{fg}_1(P\ot t^n) = \sum_{i=0} ^{n-1} g(t)^{n-i-1} P
f(t)^i.
$$
It is easy to check that $X_*(A_f^g)$ is exact in the following
cases:

\smallskip

\item{a)} $f = \al^{-r}$ with $r\ne 0$ and $g = \id$,

\smallskip

\item{b)} $f = \al^{-r-1}$ with $r\ne 0$ and $g = \al^{-1}$,

\smallskip

\item{c)} $f = \al^{1-r}$ with $r\ne 0$ and $g = \al$.

\smallskip

\ni Using this fact, Proposition~1.5 and Remark~1.2, we get that
$\HH_*(E) = \HH_*^{(0)}(E)$.

\smallskip

As usual, we let $\fu'$ denote the derivative of $u$ respect to $t$.
Moreover, given a polynomial $P\in k[t]$, we put $\T_\la(P): = P(t+\la)
- P(t)$. Let $W_{**}(E)$ be the double complex
$$
\diagram
W_{01}(E) \dto^{\phi_0}& W_{11}(E) \lto_(0.45){\de_{11}}\dto^{\phi_1} & W_{21}(E)
\lto_(0.45){\de_{21}}\dto^{\phi_2}\\
W_{00}(E) & W_{10}(E) \lto_(0.45){\de_{10}} & W_{20}(E) \lto_(0.45){\de_{20}},
\enddiagram
$$
where
$$
\alignat2
& W_{00}(E) = \bigoplus_{i\ge 0} Ax^iy^i, && \qquad W_{01}(E) =
\bigoplus_{i\ge 0} Ax^{i+1}y^ie_2,\\
& W_{20}(E) = \bigoplus_{i\ge 0} Ax^iy^{i+1}e_1,&& \qquad W_{21}(E) =
\bigoplus_{i\ge 0} Ax^iy^ie_1e_2,\\
& W_{10}(E) = W_{00}(E)\oplus W_{20}(E),&& \qquad W_{11}(E) = W_{01}(E)
\oplus W_{21}(E),
\endalignat
$$
and
$$
\align
& \de_{10}(Px^iy^i + Qx^jy^{j+1}e_1) = -\T_{\la}(Q)x^{j+1}y^{j+1}
- Q\T_{(j+1)\la}(\fu)x^jy^j,\\
& \de_{20}(Px^iy^{i+1}e_1) = \T_{\la}(P)x^{i+1}y^{i+1} + P\T_{(i+1)\la}(\fu)x^iy^i,\\
& \de_{11}(Px^{i+1}y^ie_2 + Qx^jy^je_1e_2) = -\T_{\la}(Q)x^{j+1}y^je_2
- Q\T_{j\la}(\fu)(t+\la) x^jy^{j-1}e_2,\\
& \de_{21}(Px^iy^ie_1e_2) = \T_{\la}(P)x^{i+1}y^ie_2 + P\T_{i\la}(\fu)(t+\la)
x^iy^{i-1}e_2,\\
&\phi_0(Px^{i+1}y^ie_2) = - \T_{-\la}(P)x^{i+1}y^{i+1} + P(t-\la)\T_{i\la}(\fu)x^iy^i,\\
&\phi_1(Px^{i+1}y^ie_2) =  \T_{-\la}(P)x^{i+1}y^{i+1} - P(t-\la)\T_{i\la}(\fu)x^iy^i,\\
&\phi_1(Px^iy^ie_1e_2) = \T_{-\la}(P)x^iy^{i+1}e_1 - P(t-\la)\T_{i\la}(\fu)x^{i-1}y^ie_1
- P\T_{\la}(\fu')x^iy^i,\\
&\phi_2(Px^iy^ie_1e_2) = -\T_{-\la}(P)x^iy^{i+1}e_1 + P(t-\la)\T_{i\la}(\fu)x^{i-1}y^ie_1.
\endalign
$$

\medskip

For all $0\le i\le 2$, let $\ov{\phi}_i\:H_i(W_{*1}(E))\to
H_i(W_{*0}(E))$ be the map induced by $\phi_i$. We have following
proposition:

\proclaim{Proposition 2.2.1} The Hochschild homology of $E$ is
given by
$$
\align
&\HH_0(E) = \coker(\ov{\phi}_0),\\
&\HH_1(E) = \ker(\ov{\phi}_0) \oplus \coker(\ov{\phi}_1),\\
&\HH_2(E) = \ker(\ov{\phi}_1) \oplus \coker(\ov{\phi}_2),\\
&\HH_3(E) = \ker(\ov{\phi}_2),\\
&\HH_n(E) = 0,\qquad\text{for all $n\ge 4$.}
\endalign
$$
\endproclaim

\demo{Proof} Let $\vartheta_{**}\: Y^{(0)}_{**}(E)\to W_{**}(E)$ be
the map defined by:
$$
\align
& \vartheta_{00}(Px^iy^i) = Px^iy^i,\\
& \vartheta_{10}(Px^iy^{i+1}e_1 + Qx^jy^j\ot t^n) = Px^iy^{i+1}e_1 + nt^{n-1}Qx^jy^j,\\
& \vartheta_{20}(Px^iy^{i+1}\ot t^ne_1 + Qx^jy^j\ot t^{n_1}\ot t^{n_2})
= n(t-\la)^{n-1}Px^iy^{i+1}e_1,\\
& \vartheta_{01}(Px^{i+1}y^ie_2) = Px^{i+1}y^ie_2,\\
& \vartheta_{11}(Px^iy^ie_1e_2 + Qx^{j+1}y^j\ot t^ne_2)
= Px^iy^ie_1e_2 + n(t+\la)^{n-1}Qx^jy^je_2,\\
& \vartheta_{21}(Px^iy^i\ot t^ne_1e_2 + Qx^jy^j\ot t^{n_1}\ot t^{n_2}e_2)
= nt^{n-1}Px^iy^ie_1e_2.\\
\endalign
$$
A direct computation shows that $\vartheta_{**}$ is a map of
double complexes. We assert that it is a quasi-isomorphism. To
prove this assertion, it suffices to show that $\vartheta_{*0}$ and
$\vartheta_{*1}$ are quasi-isomorphisms. It is easy to check that there
is a map of short exact sequences
$$
\spreaddiagramcolumns{-0.4pc}
\diagram
0\rto & \displaystyle{\bigoplus_{i\ge 0}}(A\ot \ov{A}^{\ot *},b_*)w_i\rto
{}\save[]+<0pc,-0.5pc>\Drop{}\dto^(0.35){\vartheta_*^{\id \id}}\restore
& Y^{(0)}_{*0}(E) \rto\dto^{\vartheta_{*0}} &
\displaystyle{\bigoplus_{i\ge 0}}(A^{\al^{-1}}_{\al^{-1}}\ot \ov{A}^{\ot *-1},b_{*-1})w_i
\rto {}\save[]+<0pc,-0.5pc>\Drop{}\dto^(0.35){\vartheta_{*-1}^{\al^{-1}\al^{-1}}}\restore & 0\\
0 \rto & \displaystyle{\bigoplus_{i\ge 0}} X_*(A)w_i \rto & W_{*0}(E) \rto &
\displaystyle{\bigoplus_{i\ge 0}}X_{*-1}(A^{\al^{-1}}_{\al^{-1}})w_i\rto & 0,
\enddiagram
$$
where $\vartheta_*^{\id \id}$ and $\vartheta_{*-1}^{\al^{-1}
\al^{-1}}$ are instances of the quasi-isomorphism
$\vartheta_*^{fg}$ introduced above. From this fact follows
immediately that $\vartheta_{*0}$ is a quasi-isomorphism. To prove
that $\vartheta_{*1}$ also is, we can proceed in a similar way.
Hence, $\HH_*(E) = \HH^{(0)}_*(E) = H_*(W_{**}(E))$. Now, the
proposition follows from the spectral sequence of a double
complex.\qed

\enddemo

Next, we use Proposition~2.2.1 to compute $\HH_*(E)$. We
consider separately the cases $\fu\in k$ and $\fu\notin k$.

\proclaim{Proposition 2.2.2} Assume that $\fu\in k$. Then
$$
\align
& \HH_0(E) = Ax^0y^0,\\
& \HH_1(E) = Ax^0y^0\oplus \bigoplus_{i\ge 0} kx^iy^{i+1}e_1,\\
& \HH_2(E) = \bigoplus_{i\ge 0} kx^iy^{i+1}e_1 \oplus
\bigoplus_{i\ge 0} kx^iy^ie_1e_2,\\
& \HH_3(E) = \bigoplus_{i\ge 0} kx^iy^ie_1e_2,\\
& \HH_n(E) = 0\qquad \text{for all $n\ge 4$.}
\endalign
$$

\endproclaim

\demo{Proof} It suffices to prove that
$$
\alignat2
&H_0(W_{*0}(E)) = Ax^0y^0,&&\qquad H_0(W_{*1}(E)) = 0,\\
&H_1(W_{*0}(E)) = Ax^0y^0 \oplus \bigoplus_{i\ge 0} kx^iy^{i+1}e_1,&&\qquad
H_1(W_{*1}(E))  = \bigoplus_{i\ge 0} kx^iy^ie_1e_2,\\
&H_2(W_{*0}(E)) = \bigoplus_{i\ge 0} kx^iy^{i+1}e_1,
&&\qquad H_2(W_{*1}(E)) = \bigoplus_{i\ge 0} kx^iy^ie_1e_2,
\endalignat
$$
and that the maps $\ov{\phi}_i$, defined above Proposition~2.2.1 are
null. We left the details to the reader.\qed

\enddemo

In the rest of this subsection we assume that $\fu\notin k$.

For all $n\ge 0$ and $-1\le j\le n-1$, let
$$
\fU_j^n(t) := (-1)^{n-j}\sum_{0\le i_1\le \cdots \le i_{n-j}\le
j+1}\fu(t+i_1\la) \cdots \fu(t+i_{n-j}\la).
$$

\proclaim{Lemma 2.2.3} Its holds that
$$
\de_{10}\biggl(x^n y^{n+1}e_1 + \sum_{j=0}^{n-1} \fU_j^n(t)x^j
y^{j+1}e_1\biggr) = (-1)^{n+1}\T_{\la}(\fu^{n+1})(t)x^0y^0.
 $$
\endproclaim

\demo{Proof} By definition
$$
\align
\de_{10}\biggl(x^n y^{n+1}e_1 + \sum_{j=0}^{n-1} \fU_j^n(t)x^j
y^{j+1}e_1\biggr)  & = - \T_{(n+1)\la}(\fu)(t)x^ny^n\\
& - \sum_{j=0}^{n-1} \T_{\la}(\fU_j^n)(t)x^{j+1}y^{j+1}\\
& - \sum_{j=0}^{n-1} \fU_j^n(t) \T_{(j+1)\la}(\fu)(t)x^jy^j.\\
\endalign
$$
Since, $\T_{\la}(\fU_{n-1}^n)(t) = -  \T_{(n+1)\la}(\fu)(t)$ and
$ \T_{\la}(\fU_j^n)(t) = - \T_{(j+2)\la}(\fu)(t)\fU_j^{n-1}(t)$
for all $0\le j< n-1$, we obtain
$$
\align
\de_{10}\biggl(x^n y^{n+1}e_1  + \sum_{j=0}^{n-1} \fU_j^n(t)&x^j
y^{j+1}e_1\biggr)   = \sum_{j=1}^{n-1} \fu(t+(j+1)\la)\bigl(\fU_{j-1}^{n-1}(t)
-\fU_j^n(t)\bigr)x^jy^j\\
& - \sum_{j=1}^{n-1} \fu(t)\bigl(\fU_j^n(t)-\fU_{j-1}^{n-1}(t)\bigr)x^jy^j
- \fU_0^n(t) \T_{\la}(\fu)(t)x^0y^0.
\endalign
$$
Since $\fU_{j-1}^{n-1}(t)-\fU_j^n(t) = \fu(t)\fU_j^{n-1}(t)$ and
$\fU_j^n(t)- \fU_{j-1}^{n-1} (t) = \fu(t+(j+1)\la)\fU_j^{n-1}(t)$,
the right side of the above equality equals to
$$
- \fU_0^n(t) \T_{\la}(\fu)(t) x^0y^0 = (-1)^{n+1} \T_{\la}(\fu^{n+1})(t)x^0y^0,
$$
as desired.\qed
\enddemo

\proclaim{Lemma 2.2.4} We have:
$$
H_0(W_{*0}(E)) = H_1(W_{*0}(E))= \frac{Ax^0y^0}{\displaystyle{\sum_{n\ge
1} k \T_{\la}(\fu^n)(t)x^0y^0}}\quad\text{and}\quad H_2(W_{*0}(E)) = 0.
$$
\endproclaim

\demo{Proof} First, we compute $H_0(W_{*0}(E))$. From the fact that
$\T_{(i+1)\la}(\fu)\ne 0$ for all $i\ge 0$ and that the map
$P\mapsto \T_{\la}(P)$ is surjective and has kernel $k$, it follows
easily that, for each $L = \sum_{i=1}^n P_ix^iy^i$, there exists a
unique element $\wt{L} = \sum_{i=1}^n \wt{P}_i x^{i-1}y^ie_1$,
such that $L - \de_{10}(\wt{L})\in Ax^0y^0$ and $t$ divides
$\wt{P}_i$, for all $1\le i\le n$. Consider the map
$\Phi\:W_{00}(E)\to Ax^0y^0$, given by $\Phi(L) = L - \de_{10}
(\wt{L})$. It is easy to check that $\Phi$ induce an isomorphism
from $H_0(W_{*0}(E)) = W_{00}(E)/\de_{10}(W_{10}(E))$ to $Ax^0y^0/(Ax^0y^0
\cap \de_{10}(W_{10}(E)))$. Hence, to finish the computation of
$H_0(W_{*0}(E))$, we only need to prove that
$$
Ax^0y^0\cap\de_{10}(W_{10}(E)) = \sum_{n\ge 1} k \T_{\la}(\fu^n)(t)x^0y^0.
$$
By Lemma~2.2.3, $\sum_{n\ge 1}k \T_{\la}(\fu^n)(t)x^0y^0 \sub
Ax^0y^0\cap\de_{10}(W_{10}(E))$. Consider an element $L =
\sum_{j=0}^m P_jx^j y^{j+1} e_1\in W_{10}(E)$, with $P_m\ne 0$. To
prove that the converse inclusion holds it suffices to show that
if $\de_{10}(L) \in Ax^0y^0$, then
$$
\de_{10}(L)\in \sum_{n\ge 1} k \T_{\la}(\fu^n)(t)x^0y^0.\tag *
$$
We proceed by induction on $m$. Since $\de_{10}(L) \in Ax^0y^0$,
we have that $P_m\in k$. So, by Lemma~2.2.3 and the inductive
hypothesis
$$
\de_{10}\biggl(L - P_m x^my^{m+1}e_1 - P_m\sum_{j=0}^{m-1}
\fU_j^m(t)x^j y^{j+1}e_1 \biggr) \in \sum_{n\ge 1} k \T_{\la}(\fu^n)(t)x^0y^0.
$$
Applying again Lemma~2.2.3, we obtain~$(*)$, which concludes the
computation of $H_0(W_{*0}(E))$. Now, we compute $\ker(\de_{10})$.
Let $L =\sum_{i=i_0}^{i_1} P_ix^iy^{i+1}e_1$ be a non-null element
of $\bigoplus_{i\ge 0} Ax^iy^{i+1}e_1$. Assume that $P_{i_0} \ne
0$. Then,
$$
\de_{10}(L) = - P_{i_0}\T_{(i_0+1)\la}(\fu)x^{i_0}y^{i_0} +
L_1x^{i_0+1}y^{i_0+1},
$$
where $L_1\in\bigoplus_{i\ge 0} Ax^iy^i$. Consequently
$\de_{10}(L)\ne 0$, since $P_{i_0}\T_{(i_0+1)\la}(\fu)\ne 0$.
Hence, $\ker(\de_{10}) = \bigoplus_{i\ge 0} Ax^iy^i$. To finish
the proof, we must check that $\de_{20}$ is injective and its
image is $\sum_{n\ge 1} k \T_{\la}(\fu^n)(t)x^0y^0$, but
theses facts follow immediately from the above computations, since
$\de_{20}(Px^iy^{i+1}e_1) = -\de_{10}(Px^iy^{i+1}e_1)$.\qed
\enddemo

Next, we compute the homology of the second row $W_{*1}(E)$ of
$W_{**}(E)$. We will need the following lemma.

\proclaim{Lemma 2.2.5} Let $n\ge 0$. It holds that
$$
\fV_n := x^n y^ne_1e_2 + \sum_{j=0}^{n-1} \fU_{j-1}^{n-1}(t+\la)x^j
y^je_1e_2
$$
belongs to the kernel of $\de_{11}$.
\endproclaim

\demo{Proof} It is similar to the proof of Lemma~2.2.3.\qed
\enddemo

\proclaim{Lemma 2.2.6} We have:
$$
H_0(W_{*1}(E)) = 0\quad\text{and}\quad H_1(W_{*1}(E)) = H_2(W_{*1}(E)) =
\bigoplus_{n\ge 0} k\fV_n.
$$
\endproclaim

\demo{Proof} First we compute $H_0(W_{*1}(E))$. From the fact that
the map $P\mapsto \T_{\la}(P)$ is surjective, it follows easily
that $\de_{11}$ also is. Hence, $H_0(W_{*1}(E)) = 0$. We assert that
$$
\ker(\de_{11}) = W_{01}(E)\oplus \bigoplus_{n\ge 0} k\fV_n.
$$
By Lemma~2.2.5, the right side of this equality is contained in
the left side. To show the converse inclusion, it suffices to
prove that if $L = \sum_{j=0}^m P_j x^j y^j e_1e_2\in
\ker(\de_{11})$, then $L\in \bigoplus_{n\ge 0} k\fV_n$. We
assume that $P_m\ne 0$, and we proceed by induction on $m$. It is
immediate that $P_m\in k$. Thus, by Lemma~2.2.5 and the
inductive hypothesis, $L - P_m\fV_m \in \bigoplus_{n\ge 0}
k\fV_n$, which implies that $L\in \bigoplus_{n\ge 0}
k\fV_n$. To finish the proof, we must check that
$\ker(\de_{21}) = \bigoplus_{n\ge 0} k\fV_n$ and that
$\Ima(\de_{21})= W_{01}(E)$, but theses facts follow immediately from
the above computations, since $\de_{21}(Px^iy^ie_1e_2) = -
\de_{11}(Px^iy^ie_1e_2)$.\qed
\enddemo

\proclaim{Lemma 2.2.7} For $n\ge 1$ let
$$
\fW_n = -\bigl(\fu'(t) - \fu'(\la)\bigr)x^{n-1}y^ne_1 - \sum_{j=1}^{n-1}
\fU_{j-1}^{n-1}(t) \bigl(\fu'(t) - \fu'(\la)\bigr) x^{j-1}y^je_1.
$$
It holds that
$$
\align
\de_{20}(\fW_n) & = -\T_{\la}(\fu')(t)x^ny^n
- \sum_{j=1}^{n-1} \fU_{j-1}^{n-1}(t+\la)\T_{\la}(\fu')(t)x^jy^j \\
& - \fU_0^{n-1}(t)\bigl(\fu'(t)-\fu'(\la)\bigr)\T_{\la}(\fu)(t)x^0y^0.
\endalign
$$
\endproclaim

\demo{Proof} To abbreviate the expressions in the proof we write
$\wt{\fu}'(t) = \fu'(t) - \fu'(\la)$ and $\wt{\fu}'(t+\la) =
\fu'(t+\la) - \fu'(\la)$. Since,
$$
\align
- \T_{\la}\left(\wt{\fu}'(t)\fU_{n-2}^{n-1}(t)\right)
& = -\T_{n\la}(\fu)(t)\fu'(\la)+ \T_{\la}(\fu')(t)\sum_{i=1}^{n-1}\fu(t+i\la)\\
&- \fu(t)\fu'(t) + \fu(t+n\la)\fu'(t+\la),
\endalign
$$
we have
$$
\align
\de_{20}(\fW_n) & = -\T_{\la}(\fu')(t)x^ny^n - \wt{\fu}'(t)
\T_{n\la}(\fu)(t)x^{n-1}y^{n-1}\\
& - \sum_{j=1}^{n-1} \fU_{j-1}^{n-1}(t+\la)\wt{\fu}'(t+\la)x^jy^j
+ \sum_{j=1}^{n-1} \fU_{j-1}^{n-1}(t)\wt{\fu}'(t)x^jy^j\\
& - \sum_{j=1}^{n-1} \fU_{j-1}^{n-1}(t) \wt{\fu}'(t)
\T_{j\la}(\fu)(t)x^{j-1}y^{j-1}\\
& = -\T_{\la}(\fu')(t)x^ny^n + \sum_{i=1}^n\fu(t+i\la)\T_{\la}(\fu')(t)x^{n-1}y^{n-1}\\
& - \sum_{j=1}^{n-2} \fU_{j-1}^{n-1}(t+\la)\wt{\fu}'(t+\la)x^jy^j
+ \sum_{j=1}^{n-2} \fU_{j-1}^{n-1}(t)\wt{\fu}'(t)x^jy^j\\
& - \sum_{j=0}^{n-2} \fU_j^{n-1}(t) \wt{\fu}'(t)\T_{(j+1)\la}(\fu)(t)x^jy^j.
\endalign
$$
Since $\fU_{j-1}^{n-1}(t)-\fu(t+(j+1)\la)\fU_j^{n-1}(t) = \fU_j^n(t)$ and
$\fU_j^n(t)+\fu(t)\fU_j^{n-1}(t) = \fU_{j-1}^{n-1}(t+\la)$, we obtain
$$
%
\align
\de_{20}(\fW_n) & = -\T_{\la}(\fu')(t)x^ny^n + \sum_{i=1}^n\fu(t+i\la)
\T_{\la}(\fu')(t)x^{n-1}y^{n-1}\\
& - \sum_{j=1}^{n-2} \fU_{j-1}^{n-1}(t+\la)\wt{\fu}'(t+\la)x^jy^j
+ \sum_{j=1}^{n-2} \fU_j^n(t)\wt{\fu}'(t)x^jy^j\\
& + \sum_{j=0}^{n-2} \fu(t)\fU_j^{n-1}(t) \wt{\fu}'(t)x^jy^j
 -  \fU_0^{n-1}(t)\wt{\fu}'(t)\fu(t+\la)x^0y^0\\
& = -\T_{\la}(\fu')(t)x^ny^n + \sum_{i=1}^n\fu(t+i\la) \T_{\la}(\fu')(t)x^{n-1}y^{n-1}\\
& - \sum_{j=1}^{n-2} \fU_{j-1}^{n-1}(t+\la)\wt{\fu}'(t+\la)x^jy^j
+ \sum_{j=1}^{n-2} \fU_{j-1}^{n-1}(t+\la)\wt{\fu}'(t)x^jy^j\\
& -  \fU_0^{n-1}(t)\wt{\fu}'(t)\T_{\la}(\fu)(t)x^0y^0\\
& = -\T_{\la}(\fu')(t)x^ny^n - \sum_{j=1}^{n-1} \fU_{j-1}^{n-1}(t+\la)
\T_{\la}(\fu')(t)x^jy^j\\
& -  \fU_0^{n-1}(t) \wt{\fu}'(t)\T_{\la}(\fu)(t)x^0y^0,
\endalign
$$
as desired.\qed
\enddemo

Let $\fV_n$ be as in Lemma~2.2.5 and let
$$
\Psi\: \bigoplus_{n\ge 0} k\fV_n \to \frac{Ax^0y^0}
{\displaystyle{\sum_{n\ge 1} k \T_{\la}(\fu^n)(t) x^0y^0}}
$$
be the map given by $\Psi(\fV_n) = (-1)^n \ov{\T_{\la}(\fu^n\fu')(t)x^0y^0}$, where
$\ov{ax^0y^0}$ denotes the class of $ax^0y^0$ in $\frac{Ax^0y^0}
{\sum_{n\ge 1} k \T_{\la}(\fu^n)(t)x^0y^0}$, for each $a\in A$.

\proclaim{Theorem 2.2.8} The Hochschild homology of $E$ is given
by:
$$
\align
& \HH_0(E) = \frac{Ax^0y^0}{\displaystyle{\sum_{n\ge 1} k \T_{\la}(\fu^n)(t)x^0y^0}},\\
& \HH_1(E) = \coker(\Psi),\\
& \HH_2(E) = \ker(\Psi),\\
& \HH_3(E) = \bigoplus_{n\ge 0} k\fV_n,\\
& \HH_n(E) = 0\qquad \text{for all $n\ge 4$.}
\endalign
$$

\endproclaim

\demo{Proof} By Lemmas~2.2.4 and 2.2.6, the maps $\ov{\phi}_0$ and
$\ov{\phi}_2$, defined above Proposition~2.2.1, are null. We
assert that $\ov{\phi}_1$ can be identified with $\Psi$. In the proof
of Lemma~2.2.4 was show that $\ker(\de_{10}) = \bigoplus_{i\ge 0}
Ax^iy^i$. Since $\phi_1(\fV_n)\in \ker(\de_{10})$, we have
$$
\phi_1(\fV_n) = -\T_{\la}(\fu')x^n y^n - \sum_{j=0}^{n-1}
\fU_{j-1}^{n-1}(t+\la)\T_{\la}(\fu')x^jy^j.
$$
Recall that, by the proof of Lemma~2.2.4, there is an isomorphism
$$
\Phi\: \ker(\de_{10})/\Ima(\de_{20})\to \frac{Ax^0y^0}
{\displaystyle{\sum_{n\ge 1} k \T_{\la}(\fu^n)(t) x^0y^0}}.
$$
This isomorphism can be described as follows: Given $L =
\sum_{i=1}^n P_ix^iy^i\in \ker(\de_{10})$, we pick up the unique
element $\wt{L} = \sum_{i=1}^n \wt{P}_i x^{i-1}y^ie_1$, such that
$L - \de_{20}(\wt{L})\in Ax^0y^0$ and $t$ divides $\wt{P}_i$ for
all $1\le i\le n$, and we put $\Phi(\ov{L}) = \ov{L -
\de_{20}(\wt{L})}$, where $\ov{L}$ denotes the class of $L$ in
$H_1(W_{*0}(E))$ and $\ov{L - \de_{20}(\wt{L})}$ denotes the class of
$L - \de_{20}(\wt{L})$ in $\frac{Ax^0y^0}{\sum_{n\ge 1}
k \T_{\la}(\fu^n)(t) x^0y^0}$. By Lemma~2.2.7, we know that if
$L = \phi_1(\fV_n)$, then $\wt{L} = \fW_n$ and
$$
\align
L - \de_{20}(\wt{L}) & = \fU_{-1}^{n-1}(t+\la)\T_{\la}(\fu')x^0y^0\\
& + \fU_0^{n-1}(t)\bigl(\fu'(t)-\fu'(\la)\bigr)\T_{\la}(\fu)(t)x^0y^0\\
& = (-1)^n\fu^n(t+\la)\bigl(\fu'(t+\la)-\fu'(\la)\bigr) -
(-1)^n\fu^n(t)\bigl(\fu'(t)-\fu'(\la)\bigr)\\
&  = (-1)^n\T_{\la}(\fu^n\fu')(t)x^0y^0
- (-1)^n \T_{\la}(\fu^n)(t)\fu'(\la)x^0y^0.
\endalign
$$
The assertion follows immediately from this fact. Now, to finish
the proof it suffices to apply Proposition~2.2.1.\qed
\enddemo

\proclaim{Corollary 2.2.9} The following facts holds:

\roster

\smallskip

\item If $u$ has degree $1$, then $\HH_2(E)$ and $\HH_3(E)$ are vector
spaces with numerable basis and $\HH_0(E)=\HH_1(E)=0$.

\smallskip

\item If $u$ has degree $2$, then $\HH_0(E)$ and $\HH_3(E)$ are vector
spaces with numerable basis and $\HH_1(E)=\HH_2(E)=0$.

\smallskip

\item If $u$ has degree greater than $2$, then $\HH_0(E)$,
$\HH_1(E)$ and $\HH_3(E)$ are vector spaces with numerable basis
and $\HH_2(E)=0$.

\smallskip

\endroster

\noindent In all the cases $\HH_n(E)=0$, for $n>3$.
\endproclaim

\demo{Proof} It follows easily from Theorem~2.2.8, using that
$\dg\bigl(\T_{\la}(\fu^n)(t)\bigr)= n\dg(\fu)- 1$ and
$\dg\bigl(\T_{\la}(\fu^n \fu')(t)\bigr)= (n+1)\dg(\fu)-2$, where
$\dg$ denotes the degree.\qed

\enddemo

The paradigm of the algebras considered in this subsection is the
enveloping algebra $\U(\fsl(2,k)$ of the simple Lie algebra
$\fsl(2,k)$, obtained taken $\la = 2$ and $\fu = -(t-1)^2/4$. In
this case applies item~(2) of Corollary~2.2.9. The result obtained
is coherent with the celebrate Whitehead's first and second
lemmas. Corollary~2.2.9 also applies to the algebras considered in
\cite{S2}.

\specialhead 2.3. Case $\bold A = \bold k[\bold t, \bold t^{-1}]$,
$\fu\in \bold A$ and $\al(\bold t) = qt$, with $q\in \bold
k\setminus\{\bold 0\}$
\endspecialhead
Note that in this case $E$ is the algebra generated over $k$ by
the variables $t,t^{-1},x,y$ and the relations $tt^{-1} =
t^{-1}t=1$, $xt = qtx$, $yt = q^{-1}ty$ and $yx = xy + \fu -
\fu(qt)$.

\smallskip

Given a polynomial $P\in k[t]$, we put
$$
\align
& \wt{T}_i^j(P) := P(q^it) - P(q^jt),\\
& T'_{q^{-1}}(P) := P(q^{-1}t) - qP(t),\\
& T'_q(P) := P(qt) - q^{-1}P(t).
\endalign
$$
For $r\in \Z$, let $\wt{W}_{**}^{(r)}(E)$ be the double complex
$$
\diagram
{\wt{W}^{(r)}_{01}(E)} \dto^{\wt{\phi}^{(r)}_0}&
{\wt{W}^{(r)}_{11}(E)} \lto_(0.45){\wt{\de}^{(r)}_{11}}
\dto^{\wt{\phi}^{(r)}_1} & {\wt{W}^{(r)}_{21}(E)}
\lto_(0.45){\wt{\de}^{(r)}_{21}} \dto^{\wt{\phi}^{(r)}_2}\\
{\wt{W}^{(r)}_{00}(E)} & {\wt{W}^{(r)}_{10}(E)}
\lto_(0.45){\wt{\de}^{(r)}_{10}} & {\wt{W}^{(r)}_{20}(E)}
\lto_(0.45){\wt{\de}^{(r)}_{20}},
\enddiagram
$$
where
$$
\alignat2
& \wt{W}^{(r)}_{00}(E) = \bigoplus_{i,j\ge 0\atop j-i = r}
Ax^iy^j, && \qquad \wt{W}^{(r)}_{01}(E) = \bigoplus_{i,j\ge 0\atop j-i = r-1}
Ax^iy^je_2,\\
& \wt{W}^{(r)}_{20}(E) = \bigoplus_{i,j\ge 0\atop j-i = r+1} Ax^iy^je_1,
&& \qquad \wt{W}^{(r)}_{21}(E) = \bigoplus_{i,j\ge 0\atop j-i = r} Ax^iy^je_1e_2,\\
& \wt{W}^{(r)}_{10}(E) = \wt{W}^{(r)}_{00}(E)
\oplus \wt{W}^{(r)}_{20}(E), && \qquad \wt{W}^{(r)}_{11}(E)
= \wt{W}^{(r)}_{01}(E)\oplus \wt{W}^{(r)}_{21}(E),
\endalignat
$$
and
$$
\align
& \wt{\de}^{(r)}_{10}(Px^iy^j + Qx^{h}y^{l}e_1) =
-\wt{T}_1^0(Q) x^{h+1}y^{l} - Q\wt{T}_{h+1}^{-r}(\fu)
x^{h}y^{l-1} + (q^{-r}-1)tPx^iy^j,\\
& \wt{\de}^{(r)}_{20}(Px^iy^je_1) = T'_q(P)x^{i+1}y^j +
q^{-1}P\wt{T}_{i+1}^{-r}(\fu)x^iy^{j-1} + (q^{-r-1}-q^{-1})tPx^iy^je_1,\\
& \wt{\de}^{(r)}_{11}(Px^iy^je_2 + Qx^{h}y^{l}e_1e_2) =
-\wt{T}_1^0(Q)x^{h+1}y^{l}e_2
- Q\wt{T}_{i+1}^{-r+1}(\fu)x^{h}y^{l-1}e_2\\
& \phantom{\wt{\de}^{(r)}_{11}(Px^iy^je_2 + Qx^{h}y^{l}e_1e_2)}+ (q^{-r+1}-q)tPx^ix^je_2,\\
& \wt{\de}^{(r)}_{21}(Px^iy^je_1e_2) = T'_q(P)x^{i+1}y^je_2 +
q^{-1}P\wt{T}_{i+1}^{-r+1}(\fu)x^iy^{j-1}e_2 + (q^{-r}-1)tPx^iy^je_1,\\
&\wt{\phi}^{(r)}_0(Px^iy^je_2) = - \wt{T}_{-1}^0(P)x^iy^{j+1}
+ P(q^{-1}t)\wt{T}_i^0(\fu)x^{i-1}y^j,\\
&\wt{\phi}^{(r)}_1(Px^iy^je_2) =  T'_{q^{-1}}(P)x^iy^{j+1}
- P(q^{-1}t)\wt{T}_i^0(\fu)x^{i-1}y^j,\\
&\wt{\phi}^{(r)}_1(Px^iy^je_1e_2) = \wt{T}_{-1}^0(P)x^iy^{j+1}e_1
- P(q^{-1}t)\wt{T}_i^0(\fu)x^{i-1}y^je_1 - P\de_{q^{-r}}(\wt{T}_1^0(\fu))x^iy^j,\\
&\wt{\phi}^{(r)}_2(Px^iy^je_1e_2) = - T'_{q^{-1}}(P)x^iy^{j+1}e_1 +
P(q^{-1}t)\wt{T}_i^0(\fu)x^{i-1}y^je_1,
\endalign
$$
where $(n)_{q^{-r}} = \frac{q^{-rn}-1}{q^{-r} - 1}$, for all $n\in
\Z$, and $\de_{q^{-r}}$ is the linear map defined by
$\de_{q^{-r}}(t^n) = (n)_{q^{-r}}t^{n-1}$.

\proclaim{Proposition 2.3.1} The map $\wt{\vartheta}^{(r)}_{**}\:
Y^{(r)}_{**}(E)\to \wt{W}^{(r)}_{**}(E)$, defined by
$$
\align
& \wt{\vartheta}^{(r)}_{00}(Px^iy^j) = Px^iy^j,\\
& \wt{\vartheta}^{(r)}_{10}(Px^iy^je_1 + Qx^{h}y^{l}\ot t^n) =
Px^iy^je_1 + (n)_{q^{-r}}
t^{n-1}Q x^{h}y^{l},\\
& \wt{\vartheta}^{(r)}_{20}(Px^iy^j\ot t^ne_1 + Qx^{h}y^{l}\ot
t^{n_1}\ot t^{n_2}) = (n)_{q^{-r}}q^{-n+1}t^{n-1}Px^iy^je_1,\\
& \wt{\vartheta}^{(r)}_{01}(Px^iy^je_2) = Px^iy^je_2,\\
&\wt{\vartheta}^{(r)}_{11}(Px^iy^je_1e_2 + Qx^{h}y^{l}\ot
t^ne_2) = Px^iy^je_1e_2 + (n)_{q^{-r}}q^{n-1}t^{n-1}Qx^{h}y^{l}e_2,\\
& \wt{\vartheta}^{(r)}_{21}(Px^iy^i\ot t^ne_1e_2 + Qx^{h}y^{l}
\ot t^{n_1}\ot t^{n_2}e_2) = (n)_{q^{-r}}t^{n-1}Px^iy^je_1e_2,\\
\endalign
$$
is a quasi-isomorphism.
\endproclaim

\demo{Proof} Mimic the proof of Proposition~2.1.2. \qed
\enddemo

\proclaim{Corollary 2.3.2} Its holds that $\HH_*(E)=
\bigoplus_{r\in \Z} H_*(\wt{W}^{(r)}_{**}(E))$.
\endproclaim

\enddocument

\newpage


Tomando $\fu = \frac{t}{1-q} - \frac{t^{-1}}{1-q^{-}}$ se obtiene el \'algebra
$E(q)$ generada por variables $t,t^{-1},x,y$ y relaciones $tt^{-1} =
t^{-1}t=1$, $xt = qtx$, $yt = q^{-1}ty$ and $yx = xy + t-t^{-1}$. Sea $X_{**}(E)$
el complejo del Teorema 4.2.1 de ``Hochschild homology of some quantum algebras''.
Hay un isomorfismo
$$
\theta_{**}\: \bigoplus_{r\in \Z} \wt{W}_{**}^{(r)}(E(q^2))\to X_{**}(E),
$$
dado por
$$
\align
&\theta_{0v}(t^nx^iy^je_2^v) = (q-q^{-1})^{i-1} x_2^n x_1^i t^jw^v,\\
&\theta_{1v}(t^nx^iy^je_2^v) = -(q-q^{-1})^{i-1} x_2^n x_1^i t^jw^v,\\
&\theta_{1v}(t^nx^iy^je_1e_2^v) = -(q-q^{-1})^i x_2^n x_1^i t^jw^v,\\
&\theta_{2v}(t^nx^iy^je_1e_2) = (q-q^{-1})^i x_2^n x_1^i t^jw^v.
\endalign
$$


\newpage


\proclaim{Notation} Given $\ba = a_1\ot\cdots\ot a_r \in \ov{A}^{\ot^r}$
and $b\in \ov{A}$, we write
$$
b\star'\ba = \sum_{l=0}^r(-1)^l a_1\ot\cdots\ot a_l\ot b \ot
\ga^{-1}(a_{l+1})\ot \cdots\ot \ga^{-1}(a_r).
$$
\endproclaim

\head 3. Hochschild cohomology \endhead

Let $M$ be an $E$-bimodule. In this section we obtain a chain
complex, giving the Hochschild cohomology of $E$ with coefficients
in $M$, which is simpler than the canonical one of Hochschild.

\proclaim{Theorem 3.1} The Hochschild cohomology $\H^*(E,M)$, of
$E$ with coefficients in $M$, is the homology of the double
complex
$$
Y^{**}(M):=\qquad \CD Y^{01}(M) @>\partial^{11}>>
Y^{11}(M)@>\partial^{21}>> Y^{21}(M)
@>\partial^{31}>> \cdots\\
@AA\varphi^0A @AA\varphi^1A @AA\varphi_2A\\
Y^{00}(M) @>\partial^{10}>> Y^{10}(M) @>\partial_{20}>>Y_{20}(M)
@>\partial^{30}>>\cdots,\\
\endCD
$$
where
$$
\align
&Y^{n0}(M) = \Hom(\ov{A}^{\ot n}\oplus \ov{A}^{\ot n-1}e_1,M),\\
&Y^{n1}(M)=\Hom(\ov{A}^{\ot n}e_2\oplus \ov{A}^{\ot n-1}e_1e_2,M)\\
&\partial^{n0}(f)(\ba) = a_1f(\ba_{2n}) + \sum_{l=1}^{n-1}(-1)^l
f(d_l(\ba))
+ (-1)^n f(\ba_{1,n-1})a_n,\\
&\partial^{n0}(f)(\ba'e_1) = a'_1f(\ba'_{2,n-1}e_1) +
\sum^{n-2}_{l=1}(-1)^l f(d_l(\ba')e_1) \\
&\phantom{\partial^{n0}(f)(\ba'e_1)} + (-1)^{n-1}
f(\ba'_{1,n-2}e_1)
\al^{-1}(a'_{n-1}) + (-1)^{n-1} \bigl(xf(\al^{-1}(\ba'))- f(\ba')x\bigr),\\
&\partial^{n1}(f)(\ba e_2) =
a_1f(\ba_{2n}e_2)+\sum_{l=1}^{n-1}(-1)^l
f(d_l(\ba)e_2) + (-1)^n f(\ba_{1,n-1}e_2)\be^{-1}(a_n),\\
&\partial^{n1}(f)(\ba' e_1e_2) = a'_1f(\ba'_{2,n-1}e_1e_2) +
\sum^{n-2}_{l=1}(-1)^l f(d_l(\ba')e_1e_2) \\
&\phantom{\partial^{n1}}+ (-1)^{n-1}
f(\ba'_{1,n-2}e_1e_2)\ga^{-1}(a'_{n-1})
+ (-1)^{n-1} \bigl(xf(\al^{-1}(\ba')e_2)- p^{-1}f(\ba'e_2)x\bigr),\\
&\varphi^n(f)(\ba e_2)= (-1)^n \bigl(yf(\be^{-1}(\ba))- f(\ba)y\bigr),\\
&\varphi^n(f)(\ba' e_1e_2) = (-1)^n
\bigl(p^{-1}yf(\be^{-1}(\ba'e_1))- f(\ba'e_1)y\bigr) +
f\bigl((p^{-1}\fu-\al(\fu))\star' \ba\bigr),
\endalign
$$
where $\ba = a_1\ot\cdots\ot a_n$ and $\ba' = a'_1\ot\cdots\ot
a'_{n-1}$.

\endproclaim

\demo{Proof} It follows by a direct computation using that
$\H^*(E,M)$ is the cohomology of $\Hom_{E^e}(Y_{**}(\wt{E}),M)$,
where $\wt{E}$ is as in the proof of Theorem~1.1.\qed
\enddemo


\newpage


\specialhead 2.3. Case $\bold A = \bold k[\bold t, \bold t^{-1}]$,
$\fu\in \bold A$ and $\al(\bold t) = qt$, with $q\in \bold
k\setminus\{\bold 0\}$ a non root of unity
\endspecialhead
Note that in this case $E$ is the algebra generated over $k$ by
the variables $t,t^{-1},x,y$ and the relations  $tt^{-1} =
t^{-1}t=1$, $xt = qtx$, $yt = q^{-1}ty$ and $yx = xy + \fu -
\fu(qt)$.

Arguing as in the previous example it is easy to check that
$\HH_*(E) = \HH_*^{(0)}(E)$.

Given a polynomial $P\in k[t]$, we put $\wt{T}_q(P) : = P(qt) -
P(t)$. Let $\wt{W}_{**}$ be the double complex
$$
\diagram {\wt{W}_{01}} \dto^{\wt{\phi}_0}& {\wt{W}_{11}}
\lto_(0.45){\wt{\de}_{11}}\dto^{\wt{\phi}_1}
& {\wt{W}_{21}} \lto_(0.45){\wt{\de}_{21}} \dto^{\wt{\phi}_2}\\
{\wt{W}_{00}} & {\wt{W}_{10}} \lto_(0.45){\wt{\de}_{10}} &
{\wt{W}_{20}} \lto_(0.45){\wt{\de}_{20}},
\enddiagram
$$
where
$$
\alignat3 & \wt{W}_{00} = \bigoplus_{i\ge 0} Ax^iy^i, && \quad
\wt{W}_{20} =
\bigoplus_{i\ge 0} Ax^iy^{i+1}e_1,&& \quad \wt{W}_{10} = \wt{W}_{00}\oplus \wt{W}_{20}\\
& \wt{W}_{01} = \bigoplus_{i\ge 0} Ax^{i+1}y^ie_2, && \quad
\wt{W}_{21} =
\bigoplus_{i\ge 0} Ax^iy^ie_1e_2,&& \quad \wt{W}_{11} = \wt{W}_{01}\oplus \wt{W}_{21}\\
\endalignat
$$
and
$$
\align & \wt{\de}_{10}(Px^iy^i + Qx^iy^{i+1}e_1) =
-\wt{T}_{q}(Q)x^{i+1}y^{i+1}
- Q\wt{T}_{q^{i+1}}(\fu)x^iy^i,\\
& \wt{\de}_{20}(Px^iy^{i+1}e_1) = \wt{T}_{q}(P)x^{i+1}y^{i+1} + P\wt{T}_{q^{i+1}}(\fu)x^iy^i,\\
& \wt{\de}_{11}(Px^{i+1}y^ie_2 + Qx^iy^ie_1e_2) =
-\wt{T}_{q}(Q)x^{i+1}y^ie_2
- Q\wt{T}_{q^i}(\fu)(qt) x^iy^{i-1}e_2,\\
& \wt{\de}_{21}(Px^iy^ie_1e_2) = \wt{T}_{q}(P)x^{i+1}y^ie_2 +
P\wt{T}_{q^i}(\fu)(qt)
x^iy^{i-1}e_2,\\
&\wt{\phi}_0(Px^{i+1}y^ie_2) = - \wt{T}_{q^{-1}}(P)x^{i+1}y^{i+1}
+ P(q^{-1}t)\wt{T}_{q^i}(\fu)x^iy^i,\\
&\wt{\phi}_1(Px^{i+1}y^ie_2) =  \wt{T}_{q^{-1}}(P)x^{i+1}y^{i+1}
- P(q^{-1}t)\wt{T}_{q^{i}}(\fu)x^iy^i,\\
&\wt{\phi}_1(Px^iy^ie_1e_2) = \wt{T}_{q^{-1}}(P)x^iy^{i+1}e_1
- P(q^{-1}t)\wt{T}_{q^i}(\fu)x^{i-1}y^ie_1 + P\wt{T}_{q^{-1}}(\fu')x^iy^i,\\
&\wt{\phi}_2(Px^iy^ie_1e_2) = -\wt{T}_{q^{-1}}(P)x^iy^{i+1}e_1 +
P(q^{-1}t)\wt{T}_{q^i}(\fu)x^{i-1}y^ie_1.
\endalign
$$

\proclaim{Proposition 2.3.1} The map $\wt{\vartheta}_{**}\:
Y^{(0)}_{**}(E)\to \wt{W}_{**}$, defined by
$$
\align
& \wt{\vartheta}_{00}(Px^iy^i) = Px^iy^i,\\
& \wt{\vartheta}_{10}(Px^iy^{i+1}e_1 + Qx^iy^i\ot t^n) = Px^iy^{i+1}e_1 + nt^{n-1}Qx^iy^i,\\
& \wt{\vartheta}_{20}(Px^iy^{i+1}\ot t^ne_1 + Qx^iy^i\ot
t^{n_1}\ot t^{n_2})
= nq^{-n+1}t^{n-1}Px^iy^{i+1}e_1,\\
& \wt{\vartheta}_{01}(Px^{i+1}y^ie_2) = Px^{i+1}y^ie_2,\\
&\wt{\vartheta}_{11}(Px^iy^ie_1e_2 + Qx^{i+1}y^i\ot t^ne_2)
= Px^iy^ie_1e_2 + nq^{n-1}t^{n-1}Qx^iy^ie_2,\\
& \wt{\vartheta}_{21}(Px^iy^i\ot t^ne_1e_2 + Qx^iy^i\ot t^{n_1}\ot
t^{n_2}e_2)
= nt^{n-1}Px^iy^ie_1e_2,\\
\endalign
$$
is a quasi-isomorphism.
\endproclaim

\demo{Proof} Mimic the proof of Proposition~2.1.2. \qed
\enddemo

\proclaim{Corollary 2.3.2} The Hochschild homology of $E$ is the
homology of $\wt{W}_{**}$.
\endproclaim

We think that this corollary can be useful to compute the
Hochschild homology for particular values of $\fu$. Unfortunately,
we did not obtain a general result for $A=k[t,t^{-1}]$, like
Theorems~2.1.1 and 2.2.7 for $A = k[t]$.

\enddocument